\documentclass{amsart}
 \usepackage{amsmath}
 \usepackage{epsfig} 
 \usepackage{graphics} 

\def\C{\mathbb C}

\def\bcases{\begin{cases}}

\def\ecases{\end{cases}}

\newcommand{\set}[1]{\left\{ #1\right\}}

\newcommand{\To}{\longrightarrow}

\newcommand{\z}{\zeta}

\newtheorem{thm}{Theorem}
\newtheorem{main theorem}{Main Theorem}
\newtheorem{corollary}{Corollary}

\newtheorem{lemma}{Lemma}
\newtheorem{prop}{Proposition}

\newtheorem{problem 1}{Problem 1}
\newtheorem{problem 2}{Problem 2}
\newtheorem{problem 3}{Problem 3}
\theoremstyle{definition}

\newtheorem{defn}{Definition}

\newcommand{\bea}{\begin{eqnarray*}}
\newcommand{\eea}{\end{eqnarray*}}

\newcommand{\be}{\begin{equation}}
\newcommand{\ee}{\end{equation}}

\begin{document}

\title[INFINITE TYPE]
{Infinite type germs of smooth hypersurfaces in $\mathbb C^n$}

\author{John Erik Forn\ae ss, Lina Lee, Yuan Zhang }
\footnote{The first author is supported by an NSF grant.
}





\today

\begin{abstract}
In this paper we discuss germs of smooth hypersurface  in $\mathbb C^n$. We  show that if a point on the boundary has infinite D'Angelo type, then there exists a formal complex curve in the hypersurface through that point.
 \end{abstract}

\maketitle

\section{Introduction}

Let $M$ be a $\mathcal C^\infty$ real  hypersurface in $\C^n$ and
$p\in M$. Then there is a  $\mathcal C^\infty$ local defining
function $r$ of $M$ near $p$, i.e., $M=\set{r=0}\cap U$, where $U$
is some neighborhood of $p$ and $dr(p) \neq 0.$ Let
$\z:(\C,0)\To(\C^n,p)$ be a germ of a nonconstant holomorphic
curve that sends $0$ to $p$. Let $\nu(\z)$ denote the degree of
vanishing at $0$ of the function $|\z_1(t)-\z_1(0)|+\cdots +
|\z_n(t)-\z_n(0)|.$ Likewise, let $\nu(\z^*r)$  denote the order
of vanishing of $r\circ \z$ at the origin.

\begin{defn}
The D'Angelo type of $M$ at $p$, $\Delta(M,p)$ is defined as
$$
\Delta(M,p):=\sup_{\z}\frac{\nu(\z^*r)}{\nu(\z)}.
$$
We say that $p$ is a point of finite type if $\Delta(M,p)<\infty$ and
infinite type otherwise. We also denote
$$\Delta(M, p, \z):=\frac{\nu(\z^*r)}{\nu(\z)}.$$
\end{defn}

In the convergent case, D'Angelo \cite{Dangelo} proved the following theorem:

\begin{thm}[D'Angelo]\label{variety-ra}
If $M$ is a germ of a real analytic hypersurface in $\C^n$ and
$p\in M$, then $\Delta(M,p)=\infty$ if and only if there exists a
germ of a complex curve $\z(t):(\C,0)\To(\C^n,p)$ lying on $M$ and
passing through $p$.
\end{thm}

In this paper we prove the formal analogue of this result. With a
formal complex curve through $p$ we mean an expression of the form
$\z(t)=(\z_1(t),\dots,\z_n(t))$ where each $\z_j$ is a formal
power series in the complex variable $t$ and $\z(0)=p.$

\begin{thm}\label{smoothdangelo}
If $M=\{r=0\}$ is a germ of a $\mathcal C^\infty$ real
hypersurface in $\C^n$ and $p\in M$, then $\Delta(M,p)=\infty$ if
and only if there exists a nonconstant formal complex curve
$\z(t):(\C,0)\To(\C^n,p)$ such that $\z(0)=p$ and $r\circ \z$
vanishes to infinite order at $0.$
\end{thm}

The proof is similar in outline to the proof by D'Angelo, but the
formal analogue is somewhat different in some places.  Denote by
${}_n\mathcal O_0$ the ring of germs of formal complex power series at 0 in
$\C^n$. A key ingredient is the following

\begin{thm}\label{nullstellensatz}
Let $I$ be an ideal in ${}_n\mathcal O_0$. If $dim({}_n\mathcal O_0/I)= \infty$,
then
 there is a formal power curve $\z$ near 0 such that $\z^*f\sim
0$ for all $f\in I$.
\end{thm}

In section 2 we collect background results. Then in Section 3 we prove the Theorem using Theorem 3, a version of the Nullstellensatz, and in the last section we prove Theorem 3. The authors are grateful to Steven Krantz for suggesting this problem.

\section{Preliminary results}

\subsection{Infinite matrices}
Let $\mathcal M$ denote the set of infinite matrices
$$\begin{bmatrix}
a_{11}& a_{12}&  ...        \\
a_{21} &   a_{22}  &   ...    \\
...& ...& ... \\
\end{bmatrix}$$

\noindent where the $\{a_{ij}\}_{1\leq i,j <\infty}$ are complex numbers.

For $A=\{a_{ij}\} \in \mathcal M$ let $A^*=\{\overline{a_{ji}}\}$
denote the adjoint matrix. We let $\mathcal M_0$ denote those
matrices $A$ for which all sums, $\sum_{j} |a_{ij}|^2 \leq 1,
\sum_{i} |a_{ij}|^2 \leq 1.$ If $A=\{a_{ij}\},B=\{b_{k\ell}\}\in
\mathcal M$ we define the matrix product $AB=\{c_{i\ell}\}$ where
$c_{i\ell}=\sum_{k} a_{ik} b_{k\ell}$ if all these sums are
convergent to a finite complex number. If $A,B\in \mathcal M_0$
then the matrix product $AB\in \mathcal M$ always exists. Also let
$\mathcal M_1\subset \mathcal M_0$ denote the unitary matrices, i.
e. those for which $AA^*=A^*A= I=:\{\delta_{ij}\}.$

\medskip

If $A=\{a_{ij}\}\in \mathcal M,A_k=\{a_{ij}^k\} \in \mathcal M_0,k=1,\dots$ and for each $i,j, a^k_{ij}\rightarrow a_{ij}, k \rightarrow \infty$ we say that $A_k$ converges weakly to $A.$ In this case it is easy to see that $A\in \mathcal M_0.$
In particular we note that if we have a sequence $U_k\in \mathcal M_0$, after taking a subsequence, we can assume that $U_k\rightarrow U\in \mathcal M_0$ and $U_k^*\rightarrow U^*\in
\mathcal M_0.$

\medskip

We let $\mathcal M_{1,k}\subset \mathcal M_1$ denote those unitary matrices $A=\{a_{ij}\}$ for which $a_{ij}=\delta_{ij}$ whenever $\max\{i,j\}>k.$ Then the matrix ${}^kA=\{a_{ij}\}_{1\leq i,j\leq k}$ is a unitary matrix on $\C^k$ and $A$ looks like the identity matrix outside the first $k \times k$ block.

\medskip

Let $\ell^2$ denote the infinite sequences of complex numbers $c=\{c_j\}_{j\geq 1}$ with norm
$\|c\|^2=\sum_j |c_j|^2.$ Then if $A\in \mathcal M_1$, then $\|Ac\|= \|c\|.$ By weak convergence
it follows that if $A\in \mathcal M_0$ is a weak limit of matrices in $\mathcal M_1$, then we still have
$\|Ac\|\leq \|c\|.$

\subsection{D'Angelo description of defining functions}

 Let $r$ be the defining function of a smooth hypersurface around $0$ in $\C^n$.
For multiindices $J=(J_1,\dots,J_n)$,$K=(K_1,\dots,K_n)$, denote $|J|:=\sum_{j=1}^n J_j$ and set
$J< K$  using some lexicographic ordering. In other words, we say
$J < K$ if and only if either $|J|<|K|$, or $|J|=|K|$ and there
exists $k\le n$ such that $J_j=K_j$ for $j\le k$ and
$J_{k+1}<K_{k+1}$. Then we can write formally

\begin{align*}
r & \sim 2\Re h+ 4\Re\sum_{J} \sum_{K\geq*J}a_{JK}z^J\overline{z}^K\\
& \sim 2\Re h+ 4\Re\sum_{J}z^J \sum_{K\geq*J}a_{JK}\overline{z}^K\\
& \sim 2\Re h+ \sum_{J}\left|z^J +\sum_{K\geq*J}\overline{a}_{JK}z^K\right|^2- \sum_{J}\left|z^J -\sum_{K\geq*J}\overline{a}_{JK}z^K\right|^2\\
& \sim \Re h+\sum_J |f_J|^2-\sum_J|g_J|^2\\
\end{align*}
where $h$ is holomorphic, $f_J=z^J +\sum_{K\geq*J}\overline{a}_{JK}z^K$ and $g_J=z^J
-\sum_{K\geq*J}\overline{a}_{JK}z^K$.

 The key property for the above decomposition is that if we truncate $r$ at any $k$ only
finitely many terms are nonzero.

\section{Formal Power Series Case}

For any given formal complex series $g$, we use $j_k(g)$ to denote
the k-jets of $g$. A slight change of D'Angelo's theorem states as
follows:

\begin{lemma}\label{jets}
Let $\z=\z_k$ be a curve passing through the origin. If
$j_{2k\nu(\z})(\z^*r)=0$, then there is an infinite unitary matrix $U_k$
such that
\begin{equation}\label{jetk}
j_{2k\nu(\z)}(\z^*h)=j_{k\nu(\z)}(\z^*(f-U_k g))=0. \end{equation}
\end{lemma}

\begin{proof}
If $j_{2k\nu(\z)}(\z^*r)=0$, then $j_{2k\nu(\z)}(\z^*h)=0$ and
$j_{2k\nu(\z)}(\|\z^*f\|^2)=j_{2k\nu(\z)}(\|\z^*g\|^2)$. Thus
$\|j_{k\nu(\z)}\z^*f\|^2=\|j_{k\nu(\z)}\z^*g\|^2$. Note the above norms make
sense since there are only finite many terms in $f$ and $g$ with
non-vanishing $k$-jets for each $k$ according to our decomposition
of $r$. For these finitely many $f_J$'s and $g_J$'s, applying
(\cite{Dangelo} Theorem 3.5), one can find an element $\tilde U$
of the group of unitary matrices such that
$$j_{k\nu(\z)}(\z^*(f-\tilde U{g}))=0. $$
Extending $\tilde U$ to an unitary operator $U$ on $l^2$ by
letting rest of the diagonal entries to be 1 and others 0, then we
have (\ref{jetk}) holds.
\end{proof}

Denote by $M_0$ the maximal ideal in ${}_nO_0$ consisting of those
which formally vanish at 0 and define $D(I):= \dim {}_n O_0/I$.

\begin{lemma}
$D(I)<\infty$ if and only if for some integer $k, M_0^k\subset I$,
where $M_0$ is a maximal ideal of ${}_n O_0$.
\end{lemma}
\begin{proof}
If $M_0^k\subset I$, then $D(I)\le D(M_0^k)<\infty$. If
$D(I)=\ell<\infty$, then $z_j^\ell\in I$ for $j=1,\dots,n$.
Therefore $M_0^\ell\subset I$.
\end{proof}

\begin{lemma}\label{inf}
Suppose that the domain $r<0$ has infinite type at $0.$ Then there
exists a matrix $U\in \mathcal M_0$ which is a limit of matrices  in $\mathcal M_1$ so that the ideal $I=I(h,
f-Ug,U^*f-g)$ has infinite dimension.
\end{lemma}

\begin{proof}
Since the domain has infinite type, there is a sequence of unitary
matrices $U_k$ as in the previous lemma. Let $U$ be any weak limit
of the $U_k.$ Suppose that $I$ has finite dimension. Then there
exists an integer $\ell$ so that $M_0^\ell\subset I.$ But then,
for large enough $k$, $M_0^\ell\subset I(h, f-U_kg,U^*_kf-g),$ a
contradiction.
\end{proof}

\begin{corollary}\label{equivalence}
Let $M=\{r=0\}$ be a smooth hypersurface in $\mathbb C^n.$ Let
$\z$ be a formal curve. Then
  $\z^*r\sim 0$ if and only if there exists an operator
$U\in \mathcal M_0$ which is a weak limit of operators in
$\mathcal M_1$ such that $\z^*f\sim 0$  for any $f\in I$ as in
Lemma \ref{inf}.
\end{corollary}

\begin{proof}
We only need to prove the sufficient direction, so consider any
formal curve $\z$ such that $\z^*h\sim \z^*(f-Ug)\sim
\z^*(U^*f-g)\sim 0$, Then $j_k(\z^*h)=0$ and
\begin{align*}
\|j_k(\z^*f)\|&=\|j_k(\z^*U g)\|=\|U(j_k(\z^* g))\| \\
&\le \|j_k(\z^*g)\| \\
&= \|j_k(\z^*U^*f)\|=\|U^*(j_k(\z^*f))\| \\
&\le \|j_k(\z^*f)\|.
\end{align*}
Therefore $\|j_k(\z^*f)\|=\|j_k(\z^*g)\|$ for any $k$. Letting $k$
goes to infinity, then $\z^*h\sim \|\z^*f\|^2 -\|\z^*g\|^2\sim 0$
and hence $\z^*r\sim 0$.
\end{proof}

To complete the proof of  Theorem 2, we only need to prove Theorem 3.

\section{Theorem 3}

In the context of formal complex power series case, many parallel
algebraic properties to those of convergent power series still
hold. To our interest, we mention that the local ring ${}_n O_0$
is Noetherian. Moreover, Weierstrass preparation Theorem and
division Theorem are both valid in the formal case. The reader may
check \cite{Lang}, \cite{Ruiz} etc for references.

For any proper ideal $I$ in ${}_nO_0$, following the idea of
Gunning(\cite{gunning}), we can similarly construct regular system
of coordinates $z_1,\ldots, z_n$ such that there exists an integer
$k$ satisfying:\\
\medskip
{\it 1. ${}_k{O_0}\cap I=0$;\\
2. ${}_{j-1}{O_0}[z_j]\cap I$ contains a Weierstrass polynomial
$p_j$ in $z_j$ for $j=k+1,\ldots,n$.}
\medskip

In addition if $I$ is prime, the theorem of the primitive element
guarantees that by making a linear change of coordinates in
$z_{k+1},\ldots,z_n$ plane, the quotient field ${}_n\tilde{M_0}$
of ${}_n{O_0}/I$ is an algebraic extension of ${}_kO_0$. i.e.
${}_n\tilde{M_0}={}_k\tilde{{O}}[\tilde z_{k+1}]$. We call the
above regular system of coordinates to be {\it strictly regular
system of coordinates}. Denote the discriminant of the unique
irreducible defining polynomial $p_{k+1}$ of $z_{k+1}$ by $D$.
Then for each coordinate function $z_j$, we can construct $q_j:=
D\cdot z_j-Q_j(z_{n+1})\in I\cap {}_kO[z_{k+1}, z_j]$ for some
$Q_j\in {}_kO[z_j], j={k+2}, \ldots, n$. We also use $p_j\in
{}_kO_0[z_j]$ to denote the defining Weierstrass polynomial of
$z_j$. The ideal generated by $p_{k+1},q_{k+2}\ldots,q_n$ is
called to be the associated ideal.

The following lemma shows that the associated ideal
$(p_{k+1},q_{k+2}\ldots,q_n)$ and the original ideal cannot differ
by very much.

\begin{lemma}\label{8}
There exists an integer $\nu$ such that $$D^{\nu} I\subseteq
(p_{k+1},q_{k+2}\ldots,q_n)\subseteq I.$$
\end{lemma}



If $I$ is a principle ideal, we have the following result.

\begin{lemma}\label{438}
Let $P(z_1,\dots,z_k,w)$ be a formal Weierstrass polynomial
$P=w^\ell+ \sum_{j<\ell}b_j(z_1,\dots,z_k)w^j$, $b_j(0)=0$.
Suppose that $D(z_1,\dots,z_k) \not \equiv 0$. Then there exists a
formal curve $\z$ in $\C^{k+1}$ such that $P\circ \z\sim 0$ .
\end{lemma}

\begin{proof}
After linear change of coordinates on $(z_1, \ldots, z_k)$, we
assume $D(z_1,0,\ldots,0) \not \equiv 0$. Write $D(z_1,0,..., 0)=
\sum_{j=s}^\infty a_jz_1^j$ where $a_s \neq 0$. Next we consider
truncations $P_r$ of $P$ where $r >> s.$ Then the discriminant of
$P_r$, $D_r$ as a polynomial in the elementary symmetric functions
can still be written as  $D_r(z_1, 0,\ldots, 0)= az_1^s + \cdots.$
For $z_1 \neq 0$ we can write $P_r(z_1,0,...,0)= \Pi_j
(w-\alpha_j^r(z_1)).$ The expression of $D_r$ leads to the
estimate $|\alpha_j^r-\alpha_i^r| \geq c|z_1|^{s}$ for some small
$c$. Let $\Delta(\alpha_j^r(z_1),\epsilon |z_1|^s)$ be the disk
with radius $\epsilon |z_1|^s$ centered at $\alpha_j^r(z_1)$ for
each $j$. Then on the boundary of these discs $|P_r| \geq
(\epsilon |z_1|^s)^\ell.$ So for all higher truncations with
$\rho>r>>s\ell$ we have as well that $|P_\rho|\geq (\epsilon
|z_1|)^\ell$ and that the zeroes of $P_\rho$ are contained in
these discs. We need to shrink $z_1$ for this to be true. We can
also make these discs much smaller. This means that we have
trapped the zero curves for $P_r$ and so we can take formal limits
and get a curve on which $P$ vanishes to infinite order.
\end{proof}

We have the following proposition:

\begin{prop}\label{prime}
Let $P\subset{}_n O_0$ be a prime ideal with $D(P)=\infty$.  Then
there exists a formal curve $\z$ such that $\z^*f\sim 0$ for all
$f\in P$.
\end{prop}

\begin{proof}



We first choose coordinates $z_1,\dots, z_n$ such that $I$ is
strictly regular. Note that since $D(P)=\infty$, $k\neq 0$.



Consider the associated ideals $(p_{k+1},q_{k+2}\ldots,q_n)$.
Notice the discriminant of $p_{k+1}$, $D$ is not identically zero.
Apply Lemma \ref{438} to $p_{k+1}$, then there exists a formal
power curve $\z'=(\z_1(t),\z_2(t),\dots,\z_{k+1}(t))$ such that
$\z'^*p_{k+1}\sim 0$. We mention that the curve can always be
chosen so that $D\circ (\z_1(t),\z_2(t),\dots,\z_{k}(t)) \ \slash
\hspace{-.3cm}\equiv 0$.

Next we need to add coordinates $z_{k+2}(t),...,z_{n}(t)$ so that
all functions in the ideal vanish. Recall $q_j=D
z_j-Q_j(z_{k+1})\in I,j>k+1$. Here $Q_j\in {}_nO[z_{k+1}].$ We
solve for $z_j=Q_j/D$ and let $\z_j(t)=z_j$ for $j=k+2,\ldots, n$.
Then $(p_{k+1},q_{k+2}\ldots,q_n)$ vanishes on
$\z(t)=(\z_1(t),\z_2(t),\dots,\z_{n}(t))$. We'll show $z_j=Q_j/D$
is one of the formal roots for the defining Weierstrass
polynomials $p_j\in {}_kO_0[z_j]$ and therefore $Q_j/D$ vanishes
at the origin. Indeed, for each $j>k+1$, denote by $n_j$ the
degree of $p_j\in {}_kO_0[z_j]$ with respect to $z_j$ and consider
$P_j=D^{n_j}p_j$. Replacing $Dz_j$ in $P_j$ by $q_j+Q_j$. We then
get $P_j=H(q_j)+G(Q_j)$ for some polynomials $H(\cdot), G(\cdot)$
with no constant terms. Since $P_j\in I$ and $p_{k+1}$ is the
defining polynomial for $z_{k+1}$, $G(Q_j)$ is divisible by
$p_{k+1}$. $P_j$ therefore identically vanishes on the curve
defined by $(z_1,\ldots, z_{k+1})=\z'(t), q_j(\z')=0$. $z_j=Q_j/D$
is then one of the formal roots for the defining Weierstrass
polynomials $p_j$.

Finally we need to show for any $f$ in the ideal $I$, $f$ vanishes
on the formal curve $\z(t)$. Indeed, by Lemma \ref{8}, there
exists some high power $\nu$, such that $D^{\nu} f\in
(p_{k+1},q_{k+2}\ldots,q_n)$. Therefore $(D^{\nu} f)\circ \z\sim
0$. Notice $D$ vanishes to finite order on the curve, thus
$\z^*(f)\sim 0$.
\end{proof}

In order to prove  Theorem 3, we also need the following
simplification lemmas.

\begin{lemma}\label{intersection}
If $I= I_1 \cap I_2$ and $D(I)=\infty$, then $D(I_j)=\infty$ for
some $j.$
\end{lemma}

\begin{proof}
Suppose that both $D(I_j)<\infty.$ Then for some integers
$k_1,k_2$ we have that $M_0^{k_j}\subset I_j; j=1,2.$ Hence
$M_0^{\max\set{k_1,k_2}}\subset I$. Hence $D(I)<\infty.$
\end{proof}



We also similarly define formal radical of an ideal $I$ as follows
\begin{defn}
$\sqrt I:=\{f: \text{there exists an integer} k \ \text{such that}
\ f^k\in I\}$.
\end{defn}

\begin{lemma}\label{radical}
If $D(P)=\infty$ then $D(\sqrt P)=\infty.$
\end{lemma}

\begin{proof}
Suppose that $D(\sqrt P)<\infty.$ Then for some $k, M_0^k\subset
\sqrt P.$ Let $f_1,...,f_s$ be generators of $M_0^k.$ Then
$f_i^{r_i}\in P$ for some $r_i$'s. If $f\in M_0^k,$ then $f=\sum
g_j f_j$, so $f^{r_1+\cdots+r_s}\in P.$ Then $M_0^K\subset P$ for
large enough $K.$ But then $D(P)<\infty.$
\end{proof}

Recall that an ideal $J$ is primary if whenever $xy\in J$ then $
x\in J$ and $y^n\in J$ for some $n$ or vice versa.

We can now prove Theorem \ref{nullstellensatz}.

Let $I$ be an ideal in ${}_n O_0$. If $D(I)= \infty$ then we want to show that
 there is a formal power curve $\z$ near 0 such that $\z^*f\sim
0$ for all $f\in I$.

\begin{proof}
By Lasker-Noether decomposition theorem, we can write
$I=P_1\cap\cdots\cap P_k$, where $P_j$'s are primary ideals. By
Lemma \ref{intersection}, we have $D(P_j)=\infty$ for some $j$.
Lemma \ref{radical} also implies $D(\sqrt{P_j})=\infty$. Notice
$\sqrt{P_j}$ is prime, Proposition \ref{prime} implies the
existence of a formal curve $\z$ such that $\z^*f\sim 0$ for all
$f\in \sqrt{P_j}$. Since$I\subset P_j\subset\sqrt{P_j}$,
$\z^*f\sim 0$ for all $f\in I$. The proof of theorem is thus
complete.
\end{proof}


\bigskip

\noindent John Erik Forn\ae ss\\
Mathematics Department\\
The University of Michigan\\
East Hall, Ann Arbor, MI 48109\\
USA\\
fornaess@umich.edu\\

\noindent Lina Lee\\
Mathematics Department\\
The University of Michigan\\
East Hall, Ann Arbor, MI 48109\\
USA\\
linalee@umich.edu\\

\noindent Yuan Zhang\\
Mathematics Department\\
Rutgers University\\
New Brunswick, NJ, 08854\\
USA\\
yuanz@math.rutgers.edu\\

\end{document}